\documentclass[11pt]{article}
\usepackage{graphics}
\usepackage{amsmath}
\usepackage{amsfonts}
\usepackage{latexsym}
\usepackage{enumerate}
\newcommand{\C}{{\bf C}}
\newcommand{\R}{{\bf R}}

\newcommand{\Z}{{\bf Z}}

\newcommand{\rtimes}{\times \kern -2pt
\vrule height 5.2pt depth 0pt width 0.4pt\;}

\begin{document}
\title{A note on compact solvmanifolds with\\
K\"ahler structures
}
\author{Keizo Hasegawa\\
Niigata University, JAPAN
}
\date{June, 2004
}
\maketitle

\baselineskip=15pt
\begin{center}
{1.  INTRODUCTION}
\end{center}

We know that the existence of K\"ahler structure on a compact complex manifold
imposes certain homological or even homotopical restrictions on its underlining
topological manifold. Hodge theory is of central importance in this line.
There have been recently certain extensions and progresses in this area of research.
Among them is the field of K\"ahler groups, in which the main subject to
study is the fundamental group of a compact K\"ahler manifold.
Once there was a conjecture that a non-abelian, finitely generated and
torsion-free nilpotent group (which is the fundamental group of a nilmanifold) can
not be a K\"ahler group, which is a generalized assertion of the result
([\ref{BG1}, \ref{H1}]) that a non-toral nilmanifold admits no K\"ahler structures.
A counter example to this conjecture was given by Campana [\ref{C}].
Later a detailed study of solvable K\"ahler groups was done by
Arapura and Nori [\ref{AN}]; they showed in particular that a solvable
K\"ahler group must be almost nilpotent, that is, it has a nilpotent subgroup
of finite index. On the other hand, the author stated in the paper [\ref{H2}]
a general conjecture on compact K\"ahlerian solvmanifolds: a compact
solvmanifold admits a K\"ahler structure if and only if it is a finite quotient of
a complex torus which is also a complex torus bundle over a complex torus;
and showed under some restriction that the conjecture is valid. 
\smallskip

In this note we will see that the above conjecture can be proved without any
restriction, based on the result (mentioned above) by Arapura and Nori [\ref{AN}],
and applying the argument being used in the proof of the main theorem on the
author's paper [\ref{H2}]. We also see that the Benson and Gordon's
conjecture on K\"ahler structures on a class of compact solvmanifolds [\ref{BG2}]
(so-called solvmanifolds of completely solvable type) can be proved as a special
case of our main result.
\smallskip

In this note we mean by a solvmanifold a compact homogeneous space of
solvable Lie group, that is, a compact differentiable manifold on which a
connected solvable Lie group $G$ acts transitively. We can assume, by taking
the universal covering group $\widetilde{G}$ of $G$, that a solvmanifold $M$ is
of the form $\widetilde{G}/D$, where $\widetilde{G}$ is a simply connected
solvable Lie group and $D$ is a closed subgroup of $\widetilde{G}$.
It should be noted that unless $M$ is a nilmanifold, a closed
subgroup $D$ may not be a discrete subgroup (a lattice) of $G$. However, it is
known [\ref{AU}] that a solvmanifold of general type has a solvmanifold
$\widetilde{G}/D$ with discrete isotropy subgroup $D$ as a finite covering. 
\smallskip

We recall some terminologies which we use in this note.
A solvmanifold $M = G/\Gamma$,
where $\Gamma$ is a discrete subgroup of a simply connected solvable Lie
group $G$, is {\em of completely solvable type} (or {\em of type (R)}),
if the adjoint representation of the Lie algebra $\mathfrak g$ of $G$ has
only real eigen values; and {\em of rigid type} in the sense of Auslander [\ref{AU}]
(or {\em of type (I)}), if the adjoint representation of $\mathfrak g$ has only
pure imaginary (including $0$) eigen values. It is clear that $M$ is both of
completely solvable and of rigid type if and only if $\mathfrak g$ is nilpotent,
that is, $M$ is a nilmanifold. We can see in the proof of main theorem that a
K\"ahlerian solvmanifold is of rigid type, and not of completely solvable type
unless it is a complex torus. This gives a proof of the Benson and Gordon's
conjecture as a special case of the theorem.
\smallskip

We state now our main theorem in the most general form:
\smallskip

\noindent {\bfseries Main Theorem.} {\em A compact solvmanifold admits a
K\"ahler structure if and only if it is a finite quotient of a complex torus which
has a structure of a complex torus bundle over a complex torus. In particular,
a compact solvmanifold of completely solvable type has a K\"ahler structure
if and only if it is a complex torus.}
\medskip

\begin{center}
{2.  PROOF OF MAIN THEOREM}
\end{center}
\medskip

Let $M$ be a compact solvmanifold of dimension $2m$ which admits a
K\"ahler structure. We first assume that $M$ is of the form $G/\Gamma$,
where $\Gamma$ is a lattice of a simply connected solvable Lie group $G$.
By the result of Arapura and Nori, we know that the fundamental group
$\Gamma$ of $M$ is almost nilpotent, that is, $\Gamma$ contains a nilpotent
subgroup $\Delta$ of finite index. Then we can find a normal nilpotent subgroup
$\Delta'$ of $\Gamma$, which is commensurable with $\Delta$, and thus defines
another lattice of $G$. Therefore $M' = G/\Delta'$ is a finite normal
covering of $M$, which is a solvmanifold with the fundamental group $\Delta'$,
a finitely generated torsion-free nilpotent group. According to the well-known
theorem of Mostow [\ref{MS}], $M'$ is diffeomorphic to a nilmanifold. Since $M'$
has a canonical K\"ahler structure induced from $M$, it follows from the result on
K\"ahlerian nilmanifolds ([\ref{BG1}, \ref{H1}]) that $M'$ must be a complex torus.
Hence $M$ is a finite quotient of complex torus. For a solvmanifold $M$ of
general type, $M$ has a solvmanifold $\widetilde{M}$ of the form $G/\Gamma$
as a finite covering, where $\Gamma$ is a lattice of a simply connected solvable
Lie group $G$. Since a solvmanifold $\widetilde{M}$ has a canonical K\"ahler
structure induced from $M$, it follows in the same way as above that $M$ must
be a finite quotient of a complex torus. 
\smallskip

We will show that $M$ has a structure of a complex torus bundle over a complex
torus. We follow the argument in the proof of the main theorem on the paper
[\ref{H2}] to which we refer for full detail.
Let $\Gamma$ be the fundamental group of $M$. Then we can express
$\Gamma$ as the extension of a finitely generated and torsion-free nilpotent
group $N$ of rank $2l$ by the free abelian group of rank $2k$, where we can
assume $2k$ is the first Betti number $b_1$ of $M$ (we have
$2k \le b_1$ in general) and $m=k+l$:
$$ 0 \rightarrow N \rightarrow \Gamma \rightarrow \Z^{2k} \rightarrow 0$$
Since $M$ is also a finite quotient of a torus $T^{2m}$,  $\Gamma$ contains
a maximal normal free abelian subgroup $\Delta$ of rank $2m$ with finite index in
$\Gamma$. We can see that $N$ must be free abelian of rank $2l$. In fact,
since $N \cap \Delta$ is a abelian subgroup of $N$ with finite index, it follows that
the real completion $\tilde{N}$ of $N$ is abelian, and thus $N$ is also abelian.
Therefore have the following
$$ 0 \rightarrow \Z^{2l} \rightarrow \Gamma \rightarrow \Z^{2k} \rightarrow 0,$$
where $\Delta = \Z^{2l} \times s_1 \Z \times s_2 \Z \times \cdots \times s_{2k} \Z$,
and $H = \Z/s_1 \Z \times \Z/s_2 \Z \times \cdots \times \Z/s_{2k} \Z$
(some of $\Z/s_i \Z$ may be trivial) is the
holonomy group of $\Gamma$. We have now that $M = \C^m/\Gamma$,
where $\Gamma$ is a Bieberbach group with holomomy group $H$. Since the action of
$H$ on $\C^l/\Z^{2l}$ is holomorphic, we see that $M$ is a holomorphic fiber bundle
over the complex torus $\C^k/\Z^{2k}$ with fiber the complex torus $\C^l/\Z^{2l}$.

Conversely, let $M$ be a finite quotient of a complex torus which is also
a complex torus bundle over a complex torus. The fundamental group $\Gamma$
of $M$ is a Bieberbach group which is expressed as the extension of a free abelian
group $\Z^{2l}$ by another free abelian group $\Z^{2k}$. As observed before,
since the action of $\Z^{2k}$ on $\Z^{2l}$ is actually the action of the finite abelian
holonomy group $H$ of $\Gamma$ on $\Z^{2l}$, and the action of $H$ on the fiber
is holomorphic, we may assume that the action of $\Z^{2k}$ is in ${\bf U}(l)$
and thus extendable to the action of $\R^{2k}$ on $\R^{2l}$,
defining a structure of solvmanifold on $M$ of the form $G/\Gamma$,
where $\Gamma$ is a lattice of a simply connected solvable Lie group $G$. In fact,
considering $\R^{2l}$ as $\C^l$, we have
$G = \C^l \rtimes \R^{2k}$ with the action
$\phi: \R^{2k} \rightarrow {\rm Aut}(\C^l)$ defined by
$$\phi(t_j) ((z_1, z_2, ..., z_l)) = 
(e^{\sqrt{-1}\,\theta^j_1 t} z_1, e^{\sqrt{-1}\,\theta^j_2 t} z_2, ...,
e^{\sqrt{-1}\, \theta^j_l t} z_l),$$
where $t_j = t e_j$ ($e_j$: the $j$-the unit vector in $\R^{2k})$, and
$e^{\sqrt{-1}\,\theta^j_i}$ is the primitive $s_j$-th root of unity,
$i = 1, 2, ...,l, j = 1, 2, ...,2k$.

As seen in the first part of the proof, we can always take $2k$ as the first Betti
number of $M$, and makes the fibration the Albanese map into the
Albanese torus. The Lie algebra $\mathfrak g$ of $G$ is the following:
$$ \mathfrak g = \{X_1, X_2, \ldots , X_{2l}, X_{2l+1}, \ldots , X_{2l+2k}\}$$
where $2k$ is the first Betti number of $M$, and the bracket multiplications are
defined by 
$$[X_{2l+i}, X_{2j-1}] = -X_{2j}, [X_{2l+i}, X_{2j}] = X_{2j-1}$$
for all $i$ with $s_i \not= 1$, $i = 1, 2, ...,2k, j = 1, 2, ..., l$,
and all other brackets vanish.
It is clear that $\mathfrak g$ is of rigid type, and is of completely solvable type
if and only if $\mathfrak g$ is abelian, that is, $M$ is a complex torus.
This completes the proof of the theorem.
\bigskip

\noindent {\bfseries Remarks.}

\begin{list}{}{\topsep=0pt \leftmargin=10pt \itemindent=5pt}

\item[1.] It is known [\ref{AA}, \ref{MR}] that a group $\Gamma$ is the
fundamental group of a flat solvmanifold, that is, a solvmanifold with flat
Riemannian metric (which may not be invariant by the Lie group action) if and
only if $\Gamma$ is an extension of a free abelian group $\Z^l$ by another
free abelian group $\Z^k$ where the action of $\Z^k$ on $\Z^l$ is finite, and
the extension corresponds to a torsion element of $H^2 (\Z^k, \Z^l)$. We can
easily check these conditions for the class of K\"ahlerain solvmanifolds we have
determined in this paper.

\item[2.] A four-dimensional solvmanifold with K\"ahler structure is nothing but
a hyperelliptic surface. It is not hard to classify hyperelliptic surfaces as
solvmanifolds [\ref{H3}].

\item[3.] According to a result of Auslander and Szczarba [\ref{AS}], a solvmanifold
$M = G/D$ has {\em the canonical torus fibration} over the torus $G/ND$ of
dimension $b_1$ (the first Betti number of $M$) with fiber a nilmanfold,
where $N$ is the nilradical of $G$ (the maximal connected normal subgroup of $G$).
Since a solvmanifold (in general) is an Eilenberg-Macline space, for a K\"ahlerian
solvmanifold $M$, we can see that the canonical torus fibration and the Albanese
map are homotopic. However it is not apriori clear that they actually coincide
(it could be a hard problem). Since that is an essential point in the proof of the
conjecture on K\"ahlerian solvmanifolds by Arapura [\ref{AR}], the proof is incomplete.
It should be also noted that by the theorem of Grauert-Fischer, a fibration (a proper
surjective holomorphic map) from a compact K\"ahler manifold to a complex torus
with fibers complex tori is locally trivial. This is well-known for K\"ahler surfaces
[\ref{BPV}], and is known to be valid in general [\ref{F}].

\item[4.] There is a recent result of Brudnyi [\ref{B}] from which we can directly see
that a K\"ahlerian solvmanifold (of general type) must be a finite
quotient of a complex torus.

\item[5.] The Benson and Gordon's conjecture may be also deduced almost directly
from the result of Apapura and Nori [\ref{AN}] and the result of Auslander [\ref{AU}]
that for a solvmanifold of the form $G/\Gamma$, where $\Gamma$ is a lattice of a
simply connected solvable Lie group $G$, if $\Gamma$ is almost nilpotent, the Lie
algebra $\mathfrak g$ is of rigid type.

\end{list}
\bigskip

\pagebreak
\begin{center}
{REFERENCES}
\end{center}
\baselineskip=15pt
\renewcommand{\labelenumi}{[\theenumi]}
\begin{enumerate}
\itemsep=0pt
\item \label{AR} D. Arapura, {\em K\"ahler solvmanifolds},
IMRN, {\bf 3} (2004), 131-137. 
\item \label{AN} D. Arapura and M. Nori, {\em Solvable fundamental groups of
algebraic varieties and K\"ahler manifolds}, Compositio Math., {\bf 116} (1999),
no. 2, 173-188. 
\item \label{AA} L. Auslander and M. Auslander, {\em Solvable Lie groups and locally
Euclidean Riemannian spaces}, Proc. Amer. Math. Soc., {\bf 9} (1958), 933-941.
\item \label{AU} L. Auslander, {\em An exposition of the structure of
solvmanifolds I II}, Bull. Amer. Math. Soc., {\bf 79} (1973),
No. 2, 227-261, 262-285.
\item \label{AS} L. Auslander and H. Szczarba, {\em Vector bundles over
noncompact solvmanifolds}, Amer. J. of Math., {\bf 97} (1975), 260-281.
\item \label{B} A. Brudnyi, {\em Solvable matrix representations of K\"ahler groups},
Differential Geometry and its Applications, {\bf 19} (2003), 167-191.
\item \label{BG1} C. Benson and C. S. Gordon, {\em K\"ahler and symplectic
structures on nilmanifolds}, Topology, {\bf 27} (1988), 755-782.
\item \label{BG2} C. Benson and C. S. Gordon, {\em K\"ahler structures on
compact solvmanifolds}, Proc. Amer. Math. Soc., {\bf 108} (1990), 971-980. 
\item \label{BPV} W. Barth, C. Peters and Van de Ven, Compact Complex Surfaces,
Ergebnisse der Mathematik und ihrer Grenzgebiete, Vol. 4, Springer-Verlag.
\item \label{C} F. Campana, {\em Remarques sur les groupes de K\"ahler nilpotents},
Ann. Sci. Ecole Norm. Sup {\bf 28} (1995), 307-316.
\item \label{F} A. Fujiki, {\em Coarse Moduli space for polirized K\"ahler manifolds},
Publ. RIMS, kyoto univ.,{\bf 20} (1984), 977-1005. 
\item \label{H1} K. Hasegawa, {\em Minimal models of nilmanifolds},
Proc. Amer. Math. Soc., {\bf 106} (1989), 65-71.
\item \label{H2} K. Hasegawa, {\em A class of compact K\"ahlerian
solvmanifolds and a general conjecture}, Geometriae Dedicata,
{\bf 78} (1999), 253-258.
\item \label{H3} K. Hasegawa, {\em Four-dimensional compact solvmanifolds
with and without complex analytic structures},
arXiv:math.CV/0401412 Jan 2004.
\item \label{MR} A. Morgan, {\em The classification of flat solvmanifolds},
Trans. Amer. Math. Soc., {\bf 239} (1978), 321-351. 
\item \label{MS} G. D. Mostow, {\em Factor spaces of solvable Lie groups},
Ann. of Math., {\bf 60} (1954), 1-27.
\end{enumerate}
\smallskip
\baselineskip12pt
\begin{flushleft}
Keizo Hasegawa\\
\vskip1pt
Department of Mathematics\\
Faculty of Education and Human Sciences\\
Niigata University, Niigata\\
JAPAN\\
\vskip2pt
E-mail: hasegawa@ed.niigata-u.ac.jp\\
\end{flushleft}

\end{document}